\RequirePackage{fix-cm}
\documentclass[smallextended]{svjour3}       
\smartqed  
\usepackage{graphicx}
\usepackage{latexsym, amsmath, amsfonts, amscd, amssymb, verbatim, longtable}
\usepackage{multirow, booktabs}

\usepackage{xcolor}

\begin{document}
	\title{Improperly Efficient Solutions in a Class of Vector Optimization Problems\thanks{This research was supported by  Vietnam Institute for Advanced Study in Mathematics (VIASM). Nguyen Thi Thu Huong was also supported by Vietnam National Foundation for Science \& and Technology Development (NAFOSTED) under grant number 06/2020/STS01.}}
	
	
	\author{N. T. T. Huong   \and N. D. Yen 
	}
	
	\authorrunning{N. T. T. Huong        \and
	N. D. Yen} 
	
	\institute{N. T. T. Huong \at
		Department of  Mathematics\\ Faculty of  Information Technology\\ Le Quy Don University\\ 100 Hoang Quoc Viet, Hanoi, Vietnam \\
		\email{nguyenhuong2308.mta@gmail.com}           
		\and
		N. D. Yen \at
		Institute of Mathematics\\ Vietnam Academy of Science and Technology\\ 18 Hoang Quoc Viet, Hanoi 10307, Vietnam\\
		\email{ndyen@math.ac.vn}}
	
	\date{Received: date / Accepted: date}

	\maketitle
	
	\begin{abstract}
		Improperly efficient solutions in the sense of Geoffrion in linear fractional vector optimization problems with unbounded constraint sets are studied in this paper. We give two sets of conditions which assure that all the efficient solutions of a given problem are improperly efficient. We also obtain necessary conditions for an efficient solution to be improperly efficient. As a result, we have new sufficient conditions for Geoffrion's proper efficiency. The obtained results enrich our knowledge on properly efficient solutions in linear fractional vector optimization.
	\end{abstract}
	
	\section{Introduction}\label{Sect_1}
	
The concept of efficiency plays a central role in vector optimization. Slightly restricted definitions of efficiency leading to proper efficiencies in various senses have  been proposed. The aim of the existing concepts of proper efficiency is to  eliminate certain efficient points that exhibit an undesirable anomaly. An efficient solution is properly efficient in the sense of Geoffrion~\cite{GE68} if there is a constant such that, for each criterion, at least one potential gain-to-loss ratio is bounded above by the constant. Latter, this concept of proper efficiency was extended by Borwein~\cite{Bo77} and Benson~\cite{Be79} to vector optimization problems, where the ordering cone  can be any nontrivial closed convex cone. When the ordering cone is the nonnegative orthant, Benson's properness is equivalent to Geoffrion's properness, while Borwein's  properness is in general weaker than Geoffrion's properness. An efficient solution that is not properly efficient is said to be improperly efficient. 

\smallskip
Studied firstly by Choo and Atkins~\cite{CA82,CA83},\textit{ linear fractional vector optimization problems} (LFVOPs) are interesting special nonconvex vector optimization problems. The importance of these problems were highlighted in~\cite[p.~203]{CA82} and \cite[Chapter~9]{S86}. Numerical methods for solving LFVOPs can be found in \cite{M95,S86}. As observed by Choo and Atkins~\cite{CA82,CA83}, the efficient solution sets of linear fractional vector optimization problems do not have the nice linear properties as in the case of linear vector optimization problems.  Based on a theorem of Robinson~\cite[Theorem 2]{Robinson_1979} on stability of monotone affine variational inequalities, several results on stability and the efficient solution sets of LFVOPs were established in~\cite{YY2011}. More information on linear fractional vector optimization problems can be found in~\cite{HYY_Optim2017,Y2012}.

\smallskip
It is well known~\cite{Choo84} that there is no difference between efficiency and Geoffrion's proper efficiency in LFVOPs with {\it bounded constraint sets}. In other words, such problems do not have any improperly efficient solution. Recently, several  results on properly efficient solutions in the sense of Geoffrion of LFVOPs with {\it unbounded constraint sets} have been obtained. In~\cite{HYY2018}, sufficient conditions for an efficient solution to be a Geoffrion's properly efficient solution are obtained by a direct approach. In \cite{HWYY2019}, other sufficient  conditions are proved by using Benson's characterization for Geoffrion's properness. In \cite{HYY2020}, sufficient conditions for an efficient solution to be a Geoffrion's properly efficient solution are established by applying some arguments of Choo~\cite{Choo84}. 

\smallskip
As shown in~\cite[Example~3.2]{HLYZ_2018}, the Borwein properly efficient solution set of a LFVOP can be strictly larger than the Geoffrion properly efficient solution set. Also, there are LFVOPs having improperly efficient solutions in the sense of Borwein (see ~\cite[Example~3.1]{HLYZ_2018}).  Verifiable sufficient conditions for an efficient point of a LFVOP to be a Borwein's properly efficient solution can be found in~\cite{HLYZ_2018}.

\smallskip
So far, the \textit{improper efficient solutions} in the sense of Geoffrion of LFVOPs with unbounded constraint sets have not been studied. In the present paper, we will obtain some sets of conditions which ensure that all the efficient solutions of a given problem are improperly efficient. Thanks to these results, some classes of abnormal LFVOPs can be described explicitly. Necessary conditions for an efficient solution to be improperly efficient will be also established. On this basis, we get new sufficient conditions for Geoffrion's proper efficiency. 

\smallskip
The outline of remaining sections is as follows. Section~\ref{Sect_2} recalls some notations, definitions, and lemmas. In Section~\ref{Sect_3}, sufficient conditions for the Geoffrion improper efficiency of a LFVOP are obtained. Section~\ref{Sect_4} establishes  necessary conditions for an efficient solution of a LFVOP to be improperly efficient in the sense of Geoffrion, and new sufficient conditions for Geoffrion's proper efficiency. Illustrative examples and comparisons of the last conditions with the preceding ones in~\cite{HYY2018} are also provided in this section. Some concluding remarks and two open questions are given in Section~\ref{Sect_5}.
	
\section{Preliminaries}\label{Sect_2}

The scalar product and the norm in the Euclidean space~$\mathbb R^n$ are denoted, respectively, by $\langle\cdot,\cdot\rangle$ and $\|\cdot\|$. Vectors in $\mathbb R^n$ are represented as rows of real numbers in the text, but they are understood as columns of real numbers in matrix calculations. If $A$ is a matrix, then $A^T$ stands for the transposed matrix. The cone generated by a subset~$D$ of an Euclidean space is denoted by ${\rm cone}\, D$, i.e., ${\rm cone}\, D=\{tx\,:\,t>0,\, x\in D\}$. The closure of ${\rm cone}\, D$ is denoted by $\overline{\rm cone}\, D$. As usual, the nonnegative orthant in $\mathbb R^m$ and the set of positive integers are denoted respectively by $\mathbb R^m_+$ and $\mathbb N$. 

\smallskip
A nonzero vector $v\in \mathbb R^n$ is said to be \cite[p.~61]{Roc70} a \textit{direction of recession} of a nonempty convex set $D\subset\mathbb R^n$ if $x+tv\in D$ for every $t\geq 0$ and every $x\in D.$ The set composed by $0\in \mathbb R^n$ and all the directions $v\in \mathbb R^n\setminus\{0\}$ satisfying the last condition, is called the \textit{recession cone} of $D$ and denoted by $0^+D.$  If $D$ is closed and convex, then $ 0^+D=\{v\in \mathbb R^n\;:\; \exists x\in D\ \, {\rm s.t.}\ \, x+tv\in D\ \, {\rm for\ all}\ \, t> 0\}.$

\begin{lemma}\label{Lemma2_v} {\rm (See \cite[Lemma~2.10]{HYY2018})} Let $D\subset\mathbb R^n$ be closed and convex, $\bar x\in D$. If $\{x^k\}$ is a sequence in $D\setminus\{\bar x\}$ with $\displaystyle\lim_{k\rightarrow\infty}\|x^k\|=+\infty$ and $\displaystyle\lim_{k\rightarrow\infty}\displaystyle\frac{x^k-\bar x}{\|x^k-\bar x\|}=v$, then $v\in 0^+D$.
\end{lemma}

Consider \textit{linear fractional functions} $f_i:\mathbb R^n\to\mathbb R,\ i=1,\dots,m$,
of the form
$$f_i(x)=\frac{a_i^Tx+\alpha_i}{b_i^Tx+\beta_i},$$ where $a_i\in \mathbb R^n, b_i\in \mathbb R^n, \alpha_i\in \mathbb R,$
and $\beta_i\in \mathbb R$. Let $K$ be a \textit{polyhedral convex set}, i.e., there exist $p\in \mathbb N$, a matrix $C\in {\mathbb R}^{p\times n},$ and a vector $d\in {\mathbb R}^p$ such that $K=\big\{x\in {\mathbb R}^n\, :\, Cx\leq d\big\}$. We assume that $K$ is nonempty and $b_i^Tx+\beta_i>0$ for all $i\in I$ and $x\in K$, where $I:=\{1,\cdots,m\}$. Put
$f(x)=(f_1(x),\dots,f_m(x))$ and let $$\Omega=\big\{x\in \mathbb R^n\,:\,b_i^Tx+\beta_i>0,\ \, \forall i\in I\big\}.$$ Clearly, $\Omega$ is open and convex, $K\subset\Omega$, and $f$ is continuously differentiable on $\Omega$. The {\it linear fractional vector optimization problem} given by $f$ and $K$ is formally written as

\smallskip
\hskip1cm (VP) \hskip 2cm ${\rm Minimize}\ \, f(x)\ \,$
subject\ to $\; x\in K$.

\begin{definition}\label{Def efficient solution} {\rm A point $x\in K$ is said to be an {\em efficient solution} (or a {\em Pareto solution}) of ${\rm (VP)}$ if $\big(f(K)-f(x)\big)\cap \big(-\mathbb R^m_+\setminus\{0\}\big)=\emptyset$.}
\end{definition}

The efficient solution set of (VP) is denoted by $E$. If $b_i=0$ and $\beta_i=1$ for all $i\in I$, then (VP) coincides with the classical multiobjective linear optimization problem.

\smallskip The next two lemmas will be used repeatedly in the sequel.

\begin{lemma}\label{Lemma1_lff} {\rm (See, e.g., \cite{M95} and \cite[Lemma~8.1]{LTY05})}
	Let $\varphi(x)=\displaystyle\frac{a^Tx+\alpha}{b^Tx+\beta}$ be a linear fractional function defined by $a,b\in \mathbb R^n$ and $\alpha,\beta\in \mathbb R$. Suppose that  $b^Tx+\beta\neq 0$ for every $x\in K_0$, where $K_0\subset\mathbb R^n$ is an arbitrary polyhedral convex set. Then, one has
	\begin{equation}\label{derivative}
	\varphi(y)-\varphi(x)=\frac{b^Tx+\beta}{b^Ty+\beta}\, \langle\nabla \varphi(x),y-x\rangle
	\end{equation} for any $x, y\in K_0$, where $\nabla \varphi(x)$
	denotes the Fr\'echet derivative of $\varphi$ at $x$.
\end{lemma}

\begin{lemma}\label{denominator-and-recession-cone} {\rm (See \cite[Lemma~3.1]{HYY2020})}
	For any $i\in I$ and $v\in 0^+K$, it holds that $b_i^Tv\geq 0$.
\end{lemma}

\begin{definition}\label{Geoffrion_efficiency}{\rm (See \cite[p.~618]{GE68})\label{Def Proper efficient solution}} {\rm One says that $\bar x\in E$ is  a \textit{Geoffrion's properly efficient solution} of (VP) if there exists a scalar $M>0$ such that, for each $i\in I$, whenever $x\in K$ and $f_i(x)<f_i(\bar x)$ one can find an index $j\in I$ such that $f_j(x)>f_j(\bar x)$ and $A_{i,j}(\bar x,x)\leq M$ with $A_{i,j}(\bar x,x):=\displaystyle\frac{f_i(\bar x)-f_i(x)}{f_j(x)-f_j(\bar x)}$.}
\end{definition}

Geoffrion's properly efficient solution set of (VP) is denoted  by $E^{Ge}$.  

\section{Sufficient Conditions}\label{Sect_3}

One may call (VP) a \textit{pathological linear fractional vector optimization problem} if all the efficient solutions are improperly efficient in the sense of Geoffrion. 

\begin{theorem}\label{Conditions for Improperly eficient_1} 
	Suppose that there exist $k\in I$ and a vector $v\in (0^+K)\setminus\{0\}$ such that $b^T_kv= 0$ and $a^T_kv< 0$. If $b^T_jv> 0$ all $j\in I\setminus\{k\}$, then any efficient solution of {\rm (VP)} is an improperly efficient in the sense of Geoffrion. 
\end{theorem}
\proof\ Let $\bar x\in E$ be given arbitrarily. To obtain a contradiction, suppose that $\bar x\in E^{Ge}$. Then there exists $M>0$ such that for each $i\in I$, if $x\in K$ and $f_{i}(x)<f_{i}(\bar x)$, then one can find $j\in I$ such that $f_{j}(x)>f_{j}(\bar x)$ and $A_{i,j}(\bar x, x)\leq M$, where the ratio  $A_{i,j}(\bar x, x)$ has been defined in Definition~\ref{Geoffrion_efficiency}. 

Since $v\in (0^+K)\setminus\{0\}$, the vector $x_t:=\bar x+tv$ belongs to $K$ for any $t>0$. By the equality $b^T_kv= 0$ one has  
\begin{equation}\label{f_ktv}
\begin{array}{rl}  f_k(x_t)=f_k(\bar x+tv)&=\displaystyle\frac{a^T_k(\bar x+tv)+\alpha_k}{b^T_k(\bar x+tv)+\beta_k}\\&=\displaystyle\frac{a^T_k\bar x+\alpha_k}{b^T_k\bar x+\beta_k}+\displaystyle\frac{a^T_kv}{b^T_k\bar x+\beta_k}t\\&=f_k(\bar x)+\displaystyle\frac{a^T_kv}{b^T_k\bar x+\beta_k}t.\end{array}\end{equation} Since $a^T_kv<0$, this implies that 
$f_k(x_t)<f_k(\bar x)$ for any $t>0$. Hence, choosing $i=k$ and recalling the choice of $M$, we can find an index $j_t\in I\setminus\{k\}$ such that $f_{j_t}(x_t)>f_{j_t}(\bar x)$ and $A_{i,j_t}(\bar x, x_t)\leq M$. The last inequality means that 
\begin{equation}\label{eq2}
f_k(\bar x)-f_k(\bar x+tv)\leq M(f_{j_t}(\bar x+tv)-f_{j_t}(\bar x)).
\end{equation}
As $j_t\in I\setminus\{k\}$ for every $t>0$, by the Dirichlet principle we can find a sequence $\{t_\ell\}_{\ell\in \mathbb N}$ of positive numbers tending to $\infty$ such that $j_{t_\ell}=j$ for some fixed index $j\in I\setminus\{k\}$.

Now, applying Lemma \ref{Lemma1_lff} to the linear fractional function $f_j$, we have
\begin{equation}\label{eq3}
0<f_j(\bar x+tv)-f_j(\bar x)=t\displaystyle\frac{b^T_j\bar x+\beta_j}{b^T_j(\bar x+tv)+\beta_j}\left\langle\nabla f_j(\bar x),v\right\rangle
\end{equation} for all $t>0$. Therefore, combining~\eqref{eq2} with~\eqref{f_ktv} and~\eqref{eq3}, we get
\begin{equation}\label{inequality}
-\displaystyle\frac{a^T_kv}{b^T_k\bar x+\beta_k}\leq M\displaystyle\frac{b^T_j\bar x+\beta_j}{b^T_j(\bar x+t_\ell v)+\beta_j}\left\langle\nabla f_j(\bar x),v\right\rangle
\end{equation}for every $\ell\in \mathbb N$. Since $b^T_jv> 0$, passing~\eqref{inequality} to the limit as $\ell\rightarrow\infty$ gives the inequality $-\displaystyle\frac{a^T_kv}{b^T_k\bar x+\beta_k}\leq 0$, which contradicts the conditions $a^T_kv<0$ and $b^T_k\bar x+\beta_k>0$. The proof is complete. $\hfill\Box$

\begin{remark} {\rm Based on Theorem~\ref{Conditions for Improperly eficient_1}, one construct infinite number of pathological LFVOPs, where any efficient solution is an improperly efficient solution in the sense of Geoffrion.}	
\end{remark}

Using Theorem~\ref{Conditions for Improperly eficient_1}, we can revisit Example~2.6 from~\cite{HYY2018} as follows.

\begin{example}\label{Example 1_a} {\rm Consider the problem (VP) with
		$$\begin{array}{rl} & K=\big\{x=(x_1,x_2)\in {\mathbb R}^2\,:\,  x_1\geq 0,\  x_2\geq 0\big\},\\
		& f_1(x)=-x_2,\ \; f_2(x)=\displaystyle\frac{x_2}{x_1+x_2+1}.\end{array}$$
		One has $E=\{(x_1,0)\,:\, x_1\geq 0\}.$ To show that all the efficient points are improperly efficient in the sense of Geoffrion by Theorem~\ref{Conditions for Improperly eficient_1}, it suffices to choose $k=1$ and $v=(0,1)$.}
\end{example}

Clearly, the above theorem can be applied only in the case where the objective function of (VP) has at most one affine component. A natural question arises: \textit{It is possible to obtain sufficient conditions for the coincidence of the set of improperly efficient solutions with the efficient solution set when {\rm (VP)} has several affine criteria, or not?} The next theorem provides an answer to this question.

\begin{theorem}\label{Conditions for Improperly eficient_2} Suppose that there exist   $i\in I$ and $v\in (0^+K)\setminus\{0\}$ such that the following conditions are satisfied:\\
	{\rm (a)} $b^T_iv= 0$ and $a^T_iv<0$,\\
	{\rm (b)} for every $j\in I\setminus\{i\}$, either $b^T_jv>0$ or $b^T_jv=0$ and $a_j^Tv\leq 0$.\\
	Then, any efficient solution of {\rm (VP)} is an improperly efficient solution in the sense of Geoffrion.
\end{theorem}
\proof\ Let $v\in (0^+K)\setminus\{0\}$ and $i\in I$ be such that the conditions~(a) and~(b) are fulfilled. Suppose that $\bar x\in E$. Thanks to the characterization of Benson for the Geoffrion properly efficient solutions (see~\cite[Theorem~3.4]{Be79}), to have $\bar x\notin E^{Ge}$ we only need to show that 
\begin{equation}\label{Be-def}
\overline{\rm cone}\,\big(f(K)+\mathbb R^m_+-f(\bar x)\big)\cap \big({-\mathbb R^m_+}\big)\ne\{0\}.
\end{equation}
Let $\{t_k\}$ be a sequence of positive numbers such that $\displaystyle\lim_{k\rightarrow\infty}t_k=+\infty$. Put $\tau_k=t_k^{-1}$ and let $u^k=0_{\mathbb{R}^m}$ for $k\in\mathbb{N}$. Define $x^k=\bar x+t_kv$ for $k\in\mathbb{N}$. By Lemma~\ref{Lemma1_lff}, for any $\ell\in I$ and $k\in\mathbb{N}$, it holds that
\begin{align*}
\begin{array}{rl} \tau_k\big(f_\ell(x^k)+u^k_\ell-f_\ell(\bar x)\big)&=\tau_k\displaystyle\frac{b_{\ell}^T\bar x+\beta_\ell}{b_\ell^Tx^k+\beta_\ell}\ \langle\nabla f_\ell(\bar x),x^k-\bar x\rangle\\&=\tau_kt_k\displaystyle\frac{b_\ell^T\bar x+\beta_\ell}{(b_\ell^T\bar x+\beta_\ell)+t_kb_\ell^Tv}\ \langle\nabla f_{\ell}(\bar x),v\rangle.\end{array}\end{align*}
So, setting $y_\ell^k=\tau_k\big(f_\ell(x^k)+u^k_\ell-f_\ell(\bar x)\big)$ and noting that $\tau_kt_k=1$, one has
\begin{equation}\label{y_kl_1}
y_\ell^k=\displaystyle\frac{b_\ell^T\bar x+\beta_\ell}{(b_\ell^T\bar x+\beta_\ell)+t_kb_\ell^Tv}\ \langle\nabla f_{\ell}(\bar x),v\rangle\quad (\forall \ell\in I,\; \forall k\in \mathbb N).
\end{equation}
The assumptions made on (VP) guarantee that  $b_\ell^T\bar x+\beta_\ell>0$ for all $\ell\in I$ and
\begin{equation}\label{derivative_f_ell}
\nabla f_{\ell}(x)=\displaystyle\frac{(b_\ell^T x+\beta_\ell)a_\ell-(a_\ell^Tx+\alpha_\ell)b_\ell}{(b_\ell^T x+\beta_\ell)^2}\quad (\forall \ell\in I,\; \forall x\in K).
\end{equation}

By~\eqref{y_kl_1},~\eqref{derivative_f_ell}, and condition (a) we have 
\begin{equation}\label{limy_kl_2}
\displaystyle\lim_{k\rightarrow\infty}y^k_i=\displaystyle\frac{a^T_iv}{b_i^T\bar x+\beta_i}<0.
\end{equation}
Now, let $j\in I\setminus\{i\}$ be given arbitrarily. From~\eqref{y_kl_1} and~\eqref{derivative_f_ell} it follows that
\begin{equation}\label{y_kj}
y_j^k=\displaystyle\frac{1}{(b_j^T\bar x+\beta_j)+t_kb_j^Tv}\ \Big\langle\displaystyle\frac{(b_j^T\bar x+\beta_j)a_j-(a_j^T\bar x+\alpha_j)b_j}{b_j^T \bar x+\beta_j},v\Big\rangle
\end{equation}
for all $k\in \mathbb N$. By condition (b) we have 
\begin{equation}\label{limy_kl}
\displaystyle\lim_{k\rightarrow\infty}y^k_j\leq 0. 
\end{equation} Indeed, if $b^T_jv>0$, passing~\eqref{y_kj} to the limit as $k\to\infty$ gives $\displaystyle\lim_{k\rightarrow\infty}y^k_j=0$. Next, suppose that $b^T_jv=0$ and $a_j^Tv\leq 0$. Then, from~\eqref{y_kj} we deduce that 
$$
\displaystyle\lim_{k\rightarrow\infty}y^k_j=\displaystyle\frac{a^T_jv}{b_j^T\bar x+\beta_j}\leq 0.
$$ Thus,~\eqref{limy_kl} is valid for every $j\in I\setminus\{i\}$. Combining this fact with~\eqref{limy_kl_2}, we see that the sequence $\{y^k\}$ with  $y^k:=(y^k_1,y^k_2,...,y^k_m)$ for $k\in \mathbb N$ has a limit $\bar y$, and $\bar y\in {-\mathbb R^m_+}\setminus\{0\}$. In addition, since $y^k=\tau_k\big(f(x^k)+u^k-f(\bar x)\big)$ belongs to ${\rm cone}\big(f(K)+\mathbb R^m_+-f(\bar x)\big)$ for all $k\in \mathbb N$, one has $\bar y\in\overline{\rm cone}\,\big(f(K)+\mathbb R^m_+-f(\bar x)\big)$. This shows that~\eqref{Be-def} is valid and completes the proof.
$\hfill\Box$

\smallskip  Note that Theorem~\ref{Conditions for Improperly eficient_2}, which was proved by using a tool from~\cite{Be79}, encompasses Theorem~\ref{Conditions for Improperly eficient_1}.

\begin{remark} {\rm In the proof of Theorem~\ref{Conditions for Improperly eficient_2}, we have shown that if there exist $i\in I$ and $v\in (0^+K)\setminus\{0\}$ satisfying (a) and (b), then~\eqref{Be-def} holds for any $\bar x\in K$.}	
\end{remark}

\begin{remark} {\rm Theorem~\ref{Conditions for Improperly eficient_2} gives us a way to  construct infinite number of pathological linear fractional vector optimization problems with more than one affine criterion, where any efficient solution is an improperly efficient in the sense of Geoffrion.}	
\end{remark}

In \cite[Example~4.7]{HYY2018}, a linear fractional vector optimization problem with two affine citeria and one fractional criterion, which has infinitely many improperly efficient in the sense of Geoffrion, was considered.  Recently, in \cite[Example~4.5]{HWYY2019}, it was proved that all the efficient solutions are improperly efficient in the sense of Geoffrion.  By Theorem~\ref{Conditions for Improperly eficient_2} we can give a short proof for the last fact.

\begin{example} {\rm (\cite[Example~4.7]{HYY2018}; see also \cite[Example~4.5]{HWYY2019}) Consider problem (VP) with $m=3,\, n=2$,
		$$\begin{array}{rl} & K=\big\{x=(x_1,x_2)\in  {\mathbb R}^2\,:\,  x_1\geq 0,\  x_2\geq 0\big\},\\
		& f_1(x)=-x_1-x_2,\quad f_2(x)=\displaystyle\frac{x_2}{x_1+x_2+1},\quad f_3(x)=x_1-x_2.\end{array}$$
		One has $E=\big\{x=(x_1,x_2)\,:\, x_1\geq 0,\; x_2\geq 0,\; x_2<x_1+1\big\}$. To show that all the efficient points are improperly efficient in the sense of Geoffrion by Theorem~\ref{Conditions for Improperly eficient_2}, it suffices to choose $i=1$ and $v=(1,1)$.} 
\end{example}

\section{Necessary Conditions}\label{Sect_4}

Necessary conditions for (VP) to have at least one improperly efficient solution in the sense of Geoffrion are given in the following theorem, whose proof relies on some results of~\cite{Be79,Choo84} and a compactification procedure. 

\begin{theorem}\label{Necessary_conditions} If {\rm (VP)} has an improperly efficient solution $\bar x$ in the sense of Geoffrion, then exists a vector $v\in (0^+K)\setminus\{0\}$ such that at least one of the following properties is valid:\\
	{\rm (c)} There is an index $i\in I$ such that $b^T_i v= 0$ and $a^T_i v\leq 0$;\\
	{\rm (d)} $\langle\nabla f_j(\bar x),v\rangle=0$ for every $j\in I$.
\end{theorem}
\proof\ Suppose that $\bar x\in E$ is a Geoffrion's improperly efficient solution. Then, by Benson's characterization for Geoffrion's efficiency,~\eqref{Be-def} holds. Fixing any nonzero vector $w\in\overline{\rm cone}\,\big(f(K)+\mathbb R^m_+-f(\bar x)\big)\cap \big({-\mathbb R^m_+}\big)$, we have $w\leq 0$ and there is a sequence $\{w^k\}\subset\overline{\rm cone}\,\big(f(K)+\mathbb R^m_+-f(\bar x)\big)\cap \big({-\mathbb R^m_+}\big)$ tending to $w$ as $k\to\infty$. For each $k$, select $x^k\in K$, $u^k\in \mathbb R^m_+$ and $\tau_k\geq 0$ such that $w^k=\tau_k\big(f(x^k)+u^k-f(\bar x)\big)$. If $\tau_k=0$ for all $k$ belonging to an infinite subset of $\mathbb N$, then there is subsequence of $\{w^k\}$ consisting of just the zero vector. This implies that $w=0$, which is impossible. So, replacing the sequence $\{\tau_k\}$ with a subsequence, we may assume that $\tau_k>0$ for all $k$.

First, let us show that the sequence $\{x^k\}$ is unbounded. On the contrary, there is $\rho>0$ such that $\|\bar x\|\leq \rho$ and $\|x^k\|\leq\rho$ for all $k\in \mathbb N$. Define $$K_\rho=\big\{x=(x_1,\dots,x_n)\in K\,:\,  -\rho\leq x_i\leq\rho\big\}$$ and observe that  $K_\rho$ is a nonempty compact polyhedral convex set. Consider the linear fractional vector optimization problem

\smallskip
\hskip1cm ${\rm (VP)}_\rho$ \hskip 2cm ${\rm Minimize}\ \, f(x)\ \,$
subject\ to $\; x\in K_\rho$. 

\smallskip
\noindent Since $\bar x$ is an efficient solution of (VP), it is an efficient solution of ${\rm (VP)}_\rho$. So, by the compactness of $K_\rho$ and the result of Choo~\cite[p.~218]{Choo84} we can assert $\bar x$ is a properly efficient solution of ${\rm (VP)}_\rho$ in the sense of Geoffrion. Hence, thanks  to Benson's characterization for Geoffrion's efficiency~\cite[Theorem~3.4]{Be79}, we have 
\begin{equation}\label{Be-def_rho}
\overline{\rm cone}\,\big(f(K_\rho)+\mathbb R^m_+-f(\bar x)\big)\cap \big({-\mathbb R^m_+}\big)\ne\{0\}.
\end{equation} On one hand, since $w^k=\tau_k\big(f(x^k)+u^k-f(\bar x)\big)\in {\rm cone}\,\big(f(K_\rho)+\mathbb R^m_+-f(\bar x)\big)$
for $k\in \mathbb N$, we can assert that $w\in\overline{\rm cone}\,\big(f(K_\rho)+\mathbb R^m_+-f(\bar x)\big)$. On the other hand,  $w\in\big({-\mathbb R^m_+}\big)\setminus\{0\}$. Clearly, the last two inclusions contradict~\eqref{Be-def_rho}. 

We have thus proved that the sequence $\{x^k\}$ is unbounded. Replacing $\{x^k\}$ with a subsequence (if necessary), we may assume that $\|x^k\|\to+\infty$ as $k\to\infty$, and $x^k\neq\bar x$ for all $k$. Without loss of generality, we can assume that the unit vectors $v^k:=\displaystyle\frac{x^k-\bar x}{\|x^k-\bar x\|}$ converge to some $v$ with $\|v\|=1$ as $k\to\infty$. By~Lemma~\ref{Lemma2_v}, one has $v\in 0^+C$. Putting $t_k=\|x^k-\bar x\|$, we get $x^k=\bar x+t_kv^k$ for $k\in \mathbb N$. Using the equality $w^k=\tau_k\big(f(x^k)+u^k-f(\bar x)\big)$ and Lemma~\ref{Lemma1_lff}, we have
$$\label{basic_relations}
	w^k_\ell =\tau_k\Big[\displaystyle\frac{b_\ell^T\bar x+\beta_\ell}{t_k^{-1}\big(b_\ell^T\bar x+\beta_\ell\big)+b_\ell^Tv^k}\langle\nabla f_\ell(\bar x),v^k\rangle+u_\ell^k\Big]\quad (\forall\ell\in I).
$$ It follows that 
$$\label{ell_ineq}
	w^k_\ell \geq\tau_k\Big[\displaystyle\frac{b_\ell^T\bar x+\beta_\ell}{t_k^{-1}\big(b_\ell^T\bar x+\beta_\ell\big)+b_\ell^Tv^k}\langle\nabla f_\ell(\bar x),v^k\rangle\Big]\quad (\forall\ell\in I).
$$
So, one has
\begin{eqnarray}\label{ell_ineq_1}
\tau_k^{-1}w^k_\ell \geq \displaystyle\frac{b_\ell^T\bar x+\beta_\ell}{t_k^{-1}\big(b_\ell^T\bar x+\beta_\ell\big)+b_\ell^Tv^k}\langle\nabla f_\ell(\bar x),v^k\rangle\quad (\forall\ell\in I).
\end{eqnarray}

According to Lemma~\ref{denominator-and-recession-cone}, we have $b_\ell^Tv\geq 0$ for all $\ell\in I$. Thus, either $b_\ell^Tv>0$ or $b_\ell^Tv=0$.

Let $I_1=\{\ell\in I\, :\, w_\ell<0\}$ and $I_2=\{\ell\in I\, :\, w_\ell=0\}$. As~$w\in\big({-\mathbb R^m_+}\big)\setminus\{0\}$, we have $I_1\cup I_2=I$ and $I_1\neq\emptyset$.

Concerning the sequence of scalars $\{\tau_k\}$, there are two possibilities: (i) $\{\tau_k\}$ is bounded; (ii) $\{\tau_k\}$ is unbounded; 

If the sequence $\{\tau_k\}$ is bounded, it must have a convergent subsequence, which is denoted again by $\{\tau_k\}$. First, consider the situation where $\displaystyle\lim_{k\to\infty}\tau_k=\bar\tau$ with $\bar\tau>0$. For each index $\ell\in I$, if  $b_\ell^Tv>0$, then passing the inequality in~\eqref{ell_ineq_1} to the limit as $k\to\infty$ gives 
$$\label{ell_ineq_2}
	\bar \tau^{-1}w_\ell \geq \displaystyle\frac{b_\ell^T\bar x+\beta_\ell}{b_\ell^Tv}\langle\nabla f_\ell(\bar x),v\rangle.
$$ This forces $\langle\nabla f_\ell(\bar x),v\rangle\leq 0$. If $b_\ell^Tv=0$, we also have $\langle\nabla f_\ell(\bar x),v\rangle\leq 0$. Indeed, if $\langle\nabla f_\ell(\bar x),v\rangle > 0$, then by letting $k\to\infty$ from~\eqref{ell_ineq_1} we get $\bar \tau^{-1}w_\ell \geq +\infty$, which is impossible. (Observe that the standing assumption
$b_i^Tx+\beta_i>0$ for all $i\in I$ and $x\in K$ implies that $$t_k^{-1}\big(b_\ell^T\bar x+\beta_\ell\big)+b_\ell^Tv^k=t_k^{-1}[b_\ell^T(\bar x+t_kv^k)+\beta_\ell]=t_k^{-1}[b_\ell^Tx^k+\beta_\ell]>0$$ for $k\in \mathbb N$.) So, we have proved that $\langle\nabla f_\ell(\bar x),v\rangle\leq 0$ for every $\ell\in I$. If there exists some $\bar\ell\in I$ with $\langle\nabla f_{\bar \ell}(\bar x),v\rangle<0$, then by using Lemma~\ref{derivative} we can show that $f_{\bar \ell}(\bar x+tv)< f_{\bar \ell}(\bar x) $ and  $f_\ell(\bar x+tv)\leq f_\ell(\bar x)$ for every $\ell\in I$, where  $t\in (0,+\infty)$ is chosen arbitrarily. Since $\bar x+tv\in K$, we get $ \bar x\notin E$, contrary to our assumption.  Therefore, we must have $\langle\nabla f_\ell(\bar x),v\rangle=0$ for every $\ell\in I$. Thus property~(d) in the formulation of our theorem is valid. Now, suppose that $\displaystyle\lim_{k\to\infty}\tau_k=0$. Select any index $\ell\in I_1$. If $b_\ell^Tv>0$, then by letting $k\to\infty$ we obtain from~\eqref{ell_ineq_1} the absurd inequality
$$\label{ell_ineq_3}
	-\infty\geq \displaystyle\frac{b^T_\ell\bar x+\beta_\ell}{b^T_\ell v}\langle\nabla f_\ell(\bar x),v\rangle.
$$ Hence, we must have $b_\ell^Tv=0$. If $\langle\nabla f_\ell(\bar x),v\rangle>0$, then the right-hand-side of the inequality in~\eqref{ell_ineq_1} tends to $+\infty$, while the left-hand-side of that inequality tends to $-\infty$ when $k\to\infty$. This is impossible. So, one has $\langle\nabla f_\ell(\bar x),v\rangle\leq 0$. Since  $$\langle\nabla f_\ell(\bar x),v\rangle=\displaystyle\frac{(b^T_\ell\bar x+\beta_\ell)a_\ell^Tv-(a^T_\ell\bar x+\alpha_\ell)b^T_\ell v}{(b^T_\ell\bar x+\beta_\ell)^2}$$ and $b^T_\ell v=0$, the last inequality implies that $a_\ell^Tv\leq 0$. For every $\ell\in I_2$, the inequality in~\eqref{ell_ineq_1} gives us nothing. Anyway, we have proved that if $\displaystyle\lim_{k\to\infty}\tau_k=0$, then property~(c) in the formulation of our theorem is valid.

If $\{\tau_k\}$ is unbounded, by considering a subsequence (if necessary), we may assume that $\displaystyle\lim_{k\to\infty}\tau_k=+\infty$. For each index $\ell\in I$, if  $b_\ell^Tv>0$, then passing the inequality in~\eqref{ell_ineq_1} to the limit as $k\to\infty$ gives 
$$\label{ell_ineq_2}
	0 \geq \displaystyle\frac{b_\ell^T\bar x+\beta_\ell}{b_\ell^Tv}\langle\nabla f_\ell(\bar x),v\rangle.
$$ This forces $\langle\nabla f_\ell(\bar x),v\rangle\leq 0$. If $b_\ell^Tv=0$, then we must have $\langle\nabla f_\ell(\bar x),v\rangle\leq 0$. Otherwise, the inequality in~\eqref{ell_ineq_1} would yield $0 \geq +\infty$, which is impossible.  Therefore, $\langle\nabla f_\ell(\bar x),v\rangle\leq 0$ for every $\ell\in I$. If there exists some $\bar\ell\in I$ with $\langle\nabla f_{\bar \ell}(\bar x),v\rangle<0$, by using Lemma~\ref{derivative} and arguing as above, we get $ \bar x\notin E$, contrary to our assumption. Thus property~(d) in the formulation of our theorem must hold.

Summing up, we have proved that, for the chosen vector $v\in (0^+K)\setminus\{0\}$, at least one of the properties~(c) and~(d) is valid.
$\hfill\Box$

\smallskip
The following corollary is immediate from Theorem~\ref{Necessary_conditions}.

\begin{corollary}\label{conditions_for_properness} If $\bar x\in E$  and there does not exist any $v\in (0^+K)\setminus\{0\}$ such that either property {\rm (c)} or property {\rm (d)} in the formulation of Theorem~\ref{Necessary_conditions} is valid, then $\bar x\in E^{Ge}$.
\end{corollary}

In comparison with the sufficient conditions for an efficient solution of {\rm (VP)} to be a properly efficient solution in the sense of Geoffrion which were given in~\cite{HWYY2019,HYY2018,HYY2020}, Corollary~\ref{conditions_for_properness} adds \textit{a quite new set of conditions}. We will present a detailed comparison of Corollary~\ref{conditions_for_properness} with the main results of~\cite{HYY2018}, omitting the analyses of the relationships between this result and the results in~\cite{HWYY2019,HYY2020}.    

\smallskip 
Fix a point $\bar x\in E$. Following~\cite{HYY2018}, we consider the next three regularity conditions:
\begin{equation}\label{c1} \begin{cases}
There\ exist\ no\ (i,j)\in I^2,\; j\neq i,\ and\ v\in (0^+K)\setminus\{0\}\ with\\
\left\langle\nabla f_i(\bar x),v\right\rangle=0 \ \, and\ \, \left\langle \nabla f_j(\bar x),v\right\rangle=0,
\end{cases}
\end{equation}
\begin{equation}\label{c2}  \begin{cases}
There\ exist\ no\ (i,j)\in I^2,\ j\neq i,\ and\ v\in (0^+K)\setminus\{0\}\ such\ that\\
b^T_iv=0,\ \; \left\langle\nabla f_i(\bar x),v\right\rangle\leq 0,\ \; \left\langle \nabla f_j(\bar x),v\right\rangle> 0,
\end{cases}
\end{equation} and 
\begin{equation}\label{c3} \begin{cases}
There\ exist\ no\ triplet\ (i,j,k)\in I^3,\ where\ i,j,k\  are\ pairwise\\ distinct, and\ v\in (0^+K)\setminus\{0\}\ with 
\left\langle\nabla f_i(\bar x),v\right\rangle<0,\ \left\langle \nabla f_j(\bar x),v\right\rangle=0,\\ \left\langle \nabla f_k(\bar x),v\right\rangle>0.
\end{cases}
\end{equation}

\smallskip
The first main result of~\cite{HYY2018} is stated as follows.

\begin{proposition}\label{Thm1_JOGO} {\rm (See \cite[Theorem~3.1]{HYY2018})} Suppose that $m=2$. If the conditions~\eqref{c1} and~\eqref{c2} are satisfied, then $\bar x\in E^{Ge}$.
\end{proposition}

For $m=2$, the regularity condition~\eqref{c1} is equivalent to saying that there does not exist any $v\in (0^+K)\setminus\{0\}$ such that property {\rm (d)} in the formulation of Theorem~\ref{Necessary_conditions} is valid. For any $v\in (0^+K)\setminus\{0\}$ and $i\in I$, if $b^T_iv=0$, then $\left\langle\nabla f_i(\bar x),v\right\rangle\leq 0$ if and only if  $a^T_i v\leq 0$ (see the proof of Theorem~\ref{Necessary_conditions}). So, for $m=2$, if there does not exist any $v\in (0^+K)\setminus\{0\}$ such that property {\rm (c)} in the formulation of Theorem~\ref{Necessary_conditions} is valid, then the regularity condition~\eqref{c2} is satisfied. Therefore, \textit{for $m=2$, the result in Corollary~\ref{conditions_for_properness}
	is weaker than the result in Proposition~\ref{Thm1_JOGO}.}

\smallskip
The second main result of~\cite{HYY2018} reads as follows.

\begin{proposition}\label{Thm2_JOGO} {\rm (See \cite[Theorem~3.2]{HYY2018})} In the case where $m\geq 3$, if the conditions~\eqref{c1}--\eqref{c3} are satisfied, then $\bar x\in E^{Ge}$.
\end{proposition}

Fix any value $m\geq 3$. Clearly, if~\eqref{c1} is satisfied, then there does not exist any $v\in (0^+K)\setminus\{0\}$ such that property {\rm (d)} in the formulation of Theorem~\ref{Necessary_conditions} is valid. Now, if there does not exist any $v\in (0^+K)\setminus\{0\}$ such that property~{\rm (c)} in the formulation of Theorem~\ref{Necessary_conditions} is valid, then the regularity condition~\eqref{c2} is satisfied. Since the regularity condition~\eqref{c3} is not required for the assertion of Corollary~\ref{conditions_for_properness}, we can conclude that the latter and Proposition~\ref{Thm2_JOGO} are \textit{incomparable results}. In fact, they are very different each from other. Note that the verification of the assumptions of Corollary~\ref{conditions_for_properness} is simpler than checking those of Proposition~\ref{Thm2_JOGO}.

\smallskip 
Let us illustrate the applicability of Corollary~\ref{conditions_for_properness} by using it to revisit some examples in~\cite{HYY2018}, which were analyzed by the results recalled in Propositions~\ref{Thm1_JOGO} and~\ref{Thm2_JOGO}.

\begin{example}\label{CA_Example}{\rm (\cite[Example 2]{CA83}; see also \cite[Example~4.1]{HYY2018}) Consider problem (VP) with
		$ K=\big\{x=(x_1,x_2)\in  {\mathbb R}^2\,:\,  x_1\geq 2,\ 0\leq x_2\leq 4\big\},$
		$f_1(x)=\displaystyle\frac{-x_1}{x_1+x_2-1}$, and $f_2(x)=\displaystyle\frac{-x_1}{x_1-x_2+3}.$ One has 
		$E=\big\{(x_1,0)\,:\, x_1\geq 2\}\cup \{(x_1,4)\,:\, x_1\geq 2\big\}.$ Since $0^+K=\{v=(v_1,0)\;:\;v_1\geq 0\}$, $b_1=(1,1)$, and $b_1=(1,-1)$,  there does not exist any $v\in (0^+K)\setminus\{0\}$ such that property {\rm (c)} is valid. As shown in~\cite{HYY2018}, if $\bar x\in E$ and if $\langle\nabla f_j(\bar x),v\rangle=0$ for all $j\in J$, then $v=0$. So, there does not exist any $v\in (0^+K)\setminus\{0\}$ such that property {\rm (d)} is valid. Hence, by Corollary~\ref{conditions_for_properness} we have $\bar x\in E^{Ge}$ for every $\bar x\in E$.}
\end{example}

\begin{example}\label{HPY_Example2}{\rm (\cite[p.~483]{HPY05a}; see also \cite[Example~4.3]{HYY2018}) Consider problem $({\rm VP})$ where $n=m=3$,
		$$\begin{array}{rl} K=\big\{x\in \mathbb R^3\,:\, &x_1+x_2-2x_3\leq  1,\ x_1-2x_2+x_3\leq  1,\\ & -2x_1+x_2+x_3\leq 1,\ x_1+x_2+x_3\geq 1 \big\},\end{array}$$
		and
		$$f_i(x)={{-x_i+\displaystyle\frac{1}{2}}\over {x_1+x_2+x_3-\displaystyle\frac{3}{4}}}\qquad (i=1,2,3).$$
		Here one has
		$$\label{E_Ew1}\begin{array}{rl}
		E= & \{(x_1,x_2,x_3)\;:\; x_1\geq  1,\ x_3=x_2=x_1-1\}\\
		& \cup \{(x_1,x_2,x_3)\;:\; x_2\geq  1,\ x_3=x_1=x_2-1\}\\
		& \cup \{(x_1,x_2,x_3)\;:\; x_3\geq  1,\ x_2=x_1=x_3-1\}\end{array}$$
		and $0^+K=\{v=(\tau,\tau,\tau)\in \mathbb R^3\, :\, \tau\geq 0\}$. Obviously, one cannot find any vector $v\in (0^+K)\setminus\{0\}$ such that property {\rm (c)} is valid. We have
		$$
		\nabla f_1(x)=\frac{1}{p(x)}\begin{pmatrix}
		-x_2-x_3+\displaystyle\frac{1}{4}, x_1-\displaystyle\frac{1}{2},x_1-\displaystyle\frac{1}{2}
		\end{pmatrix},
		$$  where $p(x):=\big(x_1+x_2+x_3-\displaystyle\frac{3}{4}\big)^2$. Select any $\bar x=(\bar x_1,\bar x_2,\bar x_3)\in E$ with $\bar x_1\geq 1$ and $\bar x_2=\bar x_3=\bar x_1-1$. Then, it holds that $\left\langle\nabla f_1(\bar x),v\right\rangle=\displaystyle\frac{5\tau}{4p(\bar x)}>0$ for any $v=(\tau,\tau,\tau)$ with $\tau>0$. Since the data of the problem under consideration is symmetric w.r.t. the variables $x_1, x_2, x_3$, this implies that, for any $\bar x\in E$,  there does not exist any $v\in (0^+K)\setminus\{0\}$ such that property {\rm (d)} is valid. Therefore, by Corollary~\ref{conditions_for_properness} we have  $E=E^{Ge}$.}
\end{example}

The number of criteria in the next linear fractional vector optimization problem can be any integer $m\geq 2$.

\begin{example}\label{HPY_Example1} {\rm (\cite[pp.~479--480]{HPY05a}; see also \cite[Example~4.3]{HYY2018}) We consider problem $({\rm VP})$ where $n=m$, $m\geq 2$,
		$$K=\Big\{x\in \mathbb R^m\,:\,x_1\geq  0,\; x_2\geq  0,\dots,\; x_m\geq   0,  \; \sum_{k=1}^m x_k\geq  1\Big\},$$
		and
		$$f_i(x)={{-x_i+\displaystyle\frac{1}{2}}\over {\displaystyle\sum_{k=1}^mx_k-\displaystyle\frac{3}{4}}} \qquad (i=1,\dots,m).$$
		Here we have
		$$\label{E_Ew2}\begin{array}{rl}
		E= &\{(x_1,0,\dots, 0)\;:\;x_1\geq  1\}\\
		& \cup \{(0,x_2,\dots,0)\;:\;x_2\geq  1\}\\
		& \dots \ \dots \ \dots\\
		& \cup \{(0,\dots,0,x_m)\;:\;x_m\geq  1\}.\end{array}$$
		Note that $b_i=(1,\dots,1)$ for  all $i\in I$, and $0^+K=\mathbb R^m_+$. So, if $b_i^Tv>0$ for any  $v\in (0^+K)\setminus\{0\}$. Hence, one cannot find any $v\in (0^+K)\setminus\{0\}$ such that property {\rm (c)} is valid. Setting  $q(x)=\left(\displaystyle\sum_{k=1}^mx_k-\displaystyle\frac{3}{4}\right)^2$, one has 
		$$
		\nabla f_i(x)=\frac{1}{q(x)}
		\begin{pmatrix}
		x_i-\displaystyle\frac{1}{2},...,-\displaystyle\sum_{k\ne i}x_k+\displaystyle\frac{1}{4},..., x_i-\displaystyle\frac{1}{2}
		\end{pmatrix}
		$$ for any $x\in K$, where the expression $-\displaystyle\sum_{k\ne i}x_k+\displaystyle\frac{1}{4}$ is the $i-$th component of $\nabla f_i(x)$. Especially, for any $\bar x\in E$, where $\bar x=(\bar x_1,0,\dots, 0)$ and $\bar x_1\geq  1$, we get	
		$$
		\nabla f_1(\bar x) =\displaystyle\frac{1}{q(\bar x)}
		\begin{pmatrix}
		\displaystyle\frac{1}{4},\bar x_1-\displaystyle\frac{1}{2},..., \bar x_1-\displaystyle\frac{1}{2}
		\end{pmatrix}.
		$$ Clearly, all the components of $\nabla f_1(\bar x)$ are positive. So, for every $v\in 0^+K\setminus\{0\}$, one has $\left\langle\nabla f_1(\bar x),v\right\rangle>0$. Since the data of the problem in question is symmetric w.r.t. the variables $x_1,\dots, x_n$, this implies that, for any $\bar x\in E$, there does not exist any $v\in (0^+K)\setminus\{0\}$ such that property {\rm (d)} is valid. Hence, by Corollary~\ref{conditions_for_properness} we can assert that $E=E^{Ge}$.} 
\end{example}
	
\section{Conclusions}\label{Sect_5}

New results on proper efficiency in the sense of Geoffrion in linear fractional vector optimization have been obtained in this paper. Namely, we have established two sets of conditions guaranteeing that all the efficient solutions of a given problem are improperly efficient. Necessary conditions for an efficient solution to be improperly efficient are also given. As a by-product, we have a quite new set of sufficient conditions for Geoffrion's proper efficiency. Our results complement the preceding ones in~\cite{HWYY2019,HYY2018,HYY2020}.

\smallskip
The following open questions seem to be interesting. Note that the second question was asked in an equivalent form in \cite[Question (Q1)]{HYY2020}.

\smallskip
\textbf{Question 1:} \textit{How to narrow the gap between the necessary conditions in Theorem~\ref{Necessary_conditions} and the sufficient conditions in Theorem~\ref{Conditions for Improperly eficient_2}?}

\smallskip
\textbf{Question 2:} \textit{Can one find any problem of the form~{\rm (VP)}, where the set of improperly efficient solutions is nonempty and it is a proper subset of $E$, or not?}
	

	%
	
\section*{Conflict of interest}

The authors declare that they have no conflict of interest.

	

\end{document}